\documentclass[11pt]{article}
\newlength{\spacing}
\newcommand{\doublespace}{\setlength{\baselineskip}{2\spacing}}

\usepackage{amsfonts}
\usepackage{amssymb}
\usepackage{graphicx}
\usepackage{newlfont}
\newcommand{\eq}[1]{\begin{equation} #1 \end{equation}}

\newtheorem{thm}{Theorem}[section]

\newtheorem{lem}[thm]{Lemma}
\newtheorem{prop}[thm]{Proposition}

\def\sec{\setcounter{equation}{0}}

\def\lam{\lambda }

\def\ka{\kappa}
\def\lab{\label }

\def\rar{\to}
\def\Rar{\Rightarrow}
\def\inft{\infty}

\def\al{\alpha}
\def\be{\beta}

\def\ep{\epsilon}
\def\sg{\sigma}

\def\cD{{\cal D}}

\def\cP{{\cal P}}

\def\cX{{\cal X}}

\def\today{\ifcase\month\or
  January\or February\or March\or April\or May\or June\or
  July\or August\or September\or October\or November\or December\fi
  \space\number\day, \number\year}

\begin{document}
\begin{titlepage}
\begin{center}
{\bf Criticality of the Exponential Rate of Decay for the
Largest Nearest Neighbor Link in Random Geometric Graphs} \\
\vspace{0.20in} by \\
\vspace{0.2in} Bhupender Gupta \\
Department of Mathematics, Indian Institute of Technology,
Kanpur 208016, India \\
\vspace{0.1in} and \\
\vspace{0.2in} {Srikanth K. Iyer \footnote{Corresponding Author.
email: skiyer@math.iisc.ernet.in}}\\
Department of Mathematics,
Indian Institute of Science, Bangalore 560012, India. \\
\end{center}
\vspace{0.2in}
%
\sloppy
\begin{center} {\bf Abstract} \end{center}

\begin{center} \parbox{4.8in}
{Let $n$ points be placed independently in $d-$dimensional space
according to the density $f(x) = A_d e^{-\lam \|x\|^\al},$ $\lam>0,$
$x \in \Re^d,$ $d \geq 2.$ Let $d_n$ be the longest edge length of
the nearest neighbor graph on these points. We show that $(\lam^{-1}\log
n)^{1-1/\al} d_n - b_n$ converges weakly to the Gumbel distribution
where $b_n \sim \frac{(d-1)}{\lam \al} \log \: \log \:n.$ We also prove the following strong law
result for the normalized nearest neighbor distance
$\tilde{d}_n := \frac{{(\lam^{-1} \log \: n)}^{1 - 1/\al} \: d_n}
{\log \log \: n}$.
\[ \frac{d-1}{\al \lam} \leq \liminf_{n \rar \infty} \tilde{d}_n \leq \limsup_{n
\rar \infty} \tilde{d}_n \leq \frac{d}{\al \lam},
\]
almost surely.
Thus, the exponential rate of decay $\al = 1$ is critical, in the
sense that for $\al > 1,$ $d_n \rar 0,$ whereas for $\al \leq 1$,
$d_n \rar \infty$ $a.s.$ as $n \rar \infty.$   } \\
\vspace{0.2in} \today
\end{center}

\vspace{0.2in}
{\sl AMS 1991 subject classifications}: \\
\hspace*{0.5in} Primary:   60D05, 60G70\\
\hspace*{0.5in} Secondary:  05C05, 90C27\\
{\sl Keywords:} Random geometric graphs, nearest neighbor graph,
Poisson point process, largest nearest neighbor link, vertex degrees.

\end{titlepage}

\doublespace

\section{Introduction and main results\lab{s1}}
\sec In this paper we prove weak and strong law results for the
largest nearest neighbor distance of points distributed according to
the probability density function
\eq{ f(x) = A_d e^{-\lam \|x\|^\al}, \qquad \lam > 0, \al > 0, \; x
\in \mathbb{R}^d, \; d \geq 2, \lab{d dim exponential density}}
where $\| \cdot \|$ is the Euclidean $(\ell_2)$ norm on $\mathbb{R}^d$ and
\eq{A_d = \frac{\al \lam^{d/\al}\Gamma(d/2+1)}{d\pi^{d/2}
\Gamma(d/\al)}. \lab{e1a}}

If $X$ has density given by (\ref{d dim exponential density}), then
$R = \| X \|$ has density,
\begin{equation}
f_{R}(r) = \frac{\al \lam^{d/\al}}{\Gamma(d/\al)} r^{d-1} e^{-\lam
r^{\al}},\qquad 0<r<\infty, \; d\geq 2. \lab{e1}
\end{equation}
The basic object of study will be the graphs $G_n$ with vertex set
$\cX_n = \{X_1,X_2,\ldots ,X_n\}$, $n=1,2,\ldots,$ where the
vertices are independently distributed according to $f.$  Edges of
$G_n$ are formed by connecting each of the vertices in $\cX_n$ to
its nearest neighbor. The longest edge of the graph $G_n$ is denoted
by $d_n$. We shall refer to $G_n$ as the nearest neighbor graph
(NNG) on $\cX_n$ and to $d_n$ as the largest nearest neighbor
distance (LNND). For any finite subset $\cX \subset \mathbb{R}^d,$ let
$G(\cX,r)$ denote the graph with vertex set $\cX$ and edges between
all pairs of vertices that are at distances less than $r.$  Thus,
$d_n$ is the minimum $r_n$ required so that the graph $G(\cX_n,r_n)$
has no isolated nodes.

The largest nearest neighbor link has been studied in the context of
computational geometry (see Dette and Henze (1989) and Steele and
Tierney (1986)) and has applications in statistics, computer
science, biology and the physical sciences. For a detailed
description of Random Geometric Graphs, their properties and
applications, we refer the reader to Penrose (2003) and references
therein.

The asymptotic distribution of $d_n$ was derived in Penrose (1997) assuming that $f$ is uniform on the unit cube.
It is shown that if the metric is assumed to be the toroidal, and if $\theta$ is the volume of the unit ball, then
$n \theta d_n^d - b_n$ converge weakly to the Gumbel distribution, where $b_n \sim \log \: n$. Penrose (1998) showed that for normally distributed points $(\al = 2)$, $\sqrt{(2\log n)} \; d_n - b_n$
converge weakly to the Gumbel distribution, where $b_n \sim (d-1) \log
\log n.$ The above result is also shown to be true for the
longest edge of the minimal spanning tree. The notation $a_n \sim b_n$ implies
that $a_n/b_n$ converges to one as
$n \rar \infty.$ Hsing and Rootzen (2005) derive the asymptotic distribution for $d_n$
in the case $d=2,$ for a large class of densities, including elliptically contoured distributions,
distributions with independent Weibull-like marginals and distributions with parallel level curves
(which includes the densities defined by (\ref{d dim exponential density})).
Appel and Russo (1997) proved strong law results for $d_n$ for graphs on
uniform points in the $d-$dimensional unit cube. Penrose (1999)
extended this to general densities having compact support $\Omega$
for which $\min_{x \in \Omega} f(x) > 0$.

Our aim in this paper is to show that when the tail of the density
decays like an exponential or slower $(\al \leq 1)$, $d_n$ diverges,
whereas for super exponential decay of the tail, $d_n \rar 0$, a.s.
as $n \rar \infty$. Properties of the one dimensional exponential
random geometric graphs have been studied in Gupta, Iyer and
Manjunath (2005). In this case, spacings
between the ordered nodes are independent and exponentially
distributed. This allows for explicit computations of many
characteristics for the graph and both strong and weak law results
can be established.

It is often easier to study the graph $G_n$ via the NNG $P_n$ on the
set $\cP_n = \{X_1,X_2,\ldots ,X_{N_n}\}$, $n \geq 1$, where
$\{N_n\}_{n \geq 1}$ is a sequence of Poisson random variables that
are independent of the sequence $\{X_n\}_{n \geq 1}$ with $E[N_n] =
n.$  $\cP_n$ is an inhomogeneous Poisson point process with
intensity function $n \: f (\cdot)$ (see Penrose (2003), Prop. 1.5).
Note that the graphs $G_n$ and $P_n$ are coupled, since the first
$\min(n,N_n)$ vertices of the two graphs are identical. We also
assume that the random variables $N_n$ are non-decreasing, so that
$\cP_1 \subset \cP_2 \subset \cP_3 \cdots.$

Let $W_n(r_n)$ (respectively $W'_n(r_n)$) be the number of vertices
of degree $0$ (isolated nodes) in $G(\cX_n,r_n)$ (respectively
$G(\cP_n,r_n)$). Let $\theta_d$ denote the volume of the
$d-$dimensional unit ball in $\mathbb{R}^d$ and let $Po(\lam)$, denote a
Poisson distribution with mean $\lam > 0.$ In what follows we will write $\log_2 n$ for $\log
\log n$ and $\log_3 n$ for $\log \log \log n$ etc.

For any $\be \in \mathbb{R},$ let $(r_n)_{n\geq 1}$ be a sequence of edge
distances that satisfies
\begin{equation}
r_n(\lam^{-1}\log \:n)^{1-1/\al} - \frac{(d-1)}{\lam \al}\log_2n +
\frac{(d-1)}{2\lam \al}\log_3n \rar \frac{\beta}{\lam \al}, \label{rn}
\end{equation}
as $n \rar \infty.$ We now state our main results.
\begin{thm} Let $(r_n)_{n\geq 1}$ satisfy (\ref{rn}) as $n \rar \infty$.
Then,
\begin{equation}
W_n(r_n) \rar Po(e^{-\beta}/ C_d) \lab{e2}
\end{equation}
in distribution, where
\eq{ C_d = \frac{\al \theta_{d-1} (d-1)!}{2} \left(\frac{d-1}{2 \pi} \right)^{\frac{d-1}{2}}. \lab{e2a}}
%
\lab{t1}
\end{thm}
An easy consequence of the above result is the following limiting distribution
for $d_n.$
%
\begin{thm}
Let $f(.)$ be the d-dimensional density defined as in (\ref{d dim
exponential density}). Let $d_n$ be the largest nearest neighbor
link of the graph $G_n$ of $n$ i.i.d. points $\cX_n =
\{X_1,X_2,\ldots ,X_n\}$ distributed according to $f.$ Then,
\begin{equation}
\lim_{n \rar \inft}P[\lam \al(\lam^{-1}\log\:n)^{1-1/\al} \; d_n -
(d-1)\log_2n + \frac{(d-1)}{2}\log_3n \leq \beta + \log(C_d)] \rar
\exp(-e^{-\beta}). \lab{e3}
\end{equation}
\lab{t2}
\end{thm}
The above result for the case $\al = 2,$ was derived in Penrose (1998). In dimension $d=2,$ Theorem~\ref{t2} follows from Theorem~7, Hsing and Rootzen (2005) (see also Example 3). Their method is based on spatial blocking and uses a locally orthogonal coordinate system with respect to the level curves. We follow the approach in Penrose (1998) and use the Chen-Stein method.

Strong law results exist in the literature only for densities that do not vanish and whose support is bounded.
Suppose $d \geq 2,$ the density $f$ is continuous, has support $\Omega,$ and that
the boundary $\partial \Omega$ is a compact $(d-1)$-dimensional $C^2$ submanifold of $\mathbb{R}^d.$
Let $f_0 > 0$ be the essential infimum of $f$ restricted to $\Omega,$ and
$f_1 = \inf_{\partial \Omega} f.$ Then (see Theorem 7.2, Penrose (2003)),
\[ \lim_{n \rar \infty} \frac{n d_n^d}{\log n} = \max\left\{\frac{c_0}{f_0}, \frac{c_1}{f_1}\right\}, \qquad \mbox{a.s.}\]
Thus, the asymptotic behavior of the LNND depends on the (reciprocal of the)
infimum of the density, since it is in the vicinity of this infimum
that points will be sparse and hence be farthest from each other. If $f_0$ or $f_1$ is zero, then the right hand side is infinite, implying that the scaling on the left is not the appropriate one. We now state a strong law result for the largest nearest neighbor distance in our case.
\begin{thm} Let $d_n$ be the LNND of the NNG $G_n$ defined on
the collection $\cX_n$ of $n$ points distributed independently and
identically according to the density $f(\cdot)$ as defined in
(\ref{d dim exponential density}). Then, almost surely, for any $d \geq 2,$
\begin{equation} \liminf_{n \rar \infty } \frac{{(\lam^{-1} \log \: n)}^{1 - 1/\al} \:
d_n}{\log_2 \: n} \geq \frac{d-1}{\al \lam}. \label{e4}
\end{equation}
\begin{equation} \limsup_{n \rar \infty } \frac{{(\lam^{-1} \log \: n)}^{1 - 1/\al} \:
d_n}{\log_2 \: n} \leq \frac{d}{\al \lam}. \label{e4a}
\end{equation}
\label{t3} \end{thm}
%
\section{Proofs and supporting results}
\sec
For any $x \in \mathbb{R}^d,$ let $B(x,r)$ denote the open ball of radius
$r$ centered at $x.$ Let
\eq{ I(x,r) := \int_{B(x,r)} f(y) \; dy. }
For $\rho > 0,$ define $I(\rho,r) = I(\rho e, r),$ where $e$ is the
$d-$dimensional unit vector $(1,0,0,\ldots, 0).$ Due to the radial
symmetry of $f,$ $I(x,r)=I(\| x \| ,r).$ The following Lemma that
provides a large $\rho$ asymptotic for $I(\rho,r)$ will be crucial
in subsequent calculations.

\begin{lem}
Let $d \geq 2,$ and $(\rho_n)_{n\geq 1}$ and $(r_n)_{n\geq 1}$ be
sequences of positive numbers satisfying $\rho_n \rar \infty,\:r_n/\rho_n
\rar 0,$ and $r_n^2 \rho_n^{\al - 2} \rar 0,$  and $r_n
\rho_n^{\al-1} \rar \infty.$ Then,
\begin{equation}
 e^{- \lam w_1(n)} \left(\Gamma\left(\frac{d+1}{2}\right)+E_n\right) H(n) \leq I(\rho_n,r_n) \leq e^{- \lam w_2(n)}
\Gamma \left(\frac{d+1}{2}\right) H(n),
\label{eqn:asy_dens_ball}
\end{equation}
where
\begin{eqnarray}
w_1(n) & = & \left\{ \begin{array}{lc}
                       \frac{\al}{2} r_n^2(\rho_n^2 - 2r_n\rho_n)^{\frac{\al}{2} - 1}, & \qquad 0 < \al \leq 2 \\
                       \frac{\al}{2} r_n^2(\rho_n^2 + 2r_n\rho_n)^{\frac{\al}{2} - 2} [ (\al - 1) \rho_n^2 + 2 r_n \rho_n]
, & \qquad  \al > 2,
                       \end{array} \right.
\label{eqn:w_1} \\
%
w_2(n) & = & \left\{ \begin{array}{lc}
                        \frac{\al(\al-2)}{2}(r_n \rho_n)^2(\rho_n - 2r_n\rho_n)^{\frac{\al}{2} - 2}, & 0 < \al \leq 2 \\
                       0 , & \al > 2,
                       \end{array} \right.
\label{eqn:w_2} \\
%
\mid E_n \mid & \leq &  \frac{C_1}{r_n \rho_n^{\al-1}}, \\
\label{eqn:bound_E_n}
%
H(n) & = & A_d \theta_{d-1}2^{\frac{d-1}{2}}
r_n^d
\exp(-\lam(\rho_n^{\al} - \al r_n\rho_n^{\al-1}))(\lam \al r_n
\rho_n^{\al-1})^{-\frac{d+1}{2}},
\label{eqn:H_n}
\end{eqnarray}
%
%
where $A_d$ is as defined in (\ref{e1a}), $\theta_{d-1}$ is the volume
of the $(d-1)$-dimensional unit ball, and $C_1$ is some constant. As $n\rar \infty, \: E_n \rar
0,$ and  $w_i(n)\rar 0, \: i=1,2.$
\label{density over the ball}
\end{lem}
{\textbf{ Proof.}} In the definition of $I(\rho_n,r_n)=I(\rho_n e,r_n),$ set $y=
(\rho_n+r_n t, r_n s), \: t\in(-1,1),\:s \in \mathbb{R}^{d-1}.$ This gives,
\begin{equation}
I(\rho_n,r_n) = A_d \int_{-1}^{1}\int_{\tiny{\|s\|^2
\leq(1-t^2),s\in {\mathbb{R}}^{d-1}}}
\exp\left(-\lam((\rho_n+r_nt)^2 + (\|s\|r_n)^2
)^{\frac{\al}{2}}\right)r_n^d \ ds \ dt.\label{In}
\end{equation}
%
%
Consider first the case $0 < \al \leq 2.$ Using the Taylor's expansion
we get,
\begin{eqnarray}
((\rho_n+r_nt)^2 + (\|s\|r_n)^2 )^{\frac{\al}{2}}
& = &((\rho_n^2+2r_nt \rho_n) + (t^2+\|s\|^2)r_n^2 )^{\frac{\al}{2}} \nonumber \\
& = & (\rho_n^2+ 2r_n\rho_n t)^{\frac{\al}{2}} + (r_n^2(t^2 + ||s||^2))
\frac{\al}{2} (\rho_n^2+ 2r_n\rho_n t+\xi)^{\frac{\al}{2} -1} \nonumber \\
& = & (\rho_n^2+ 2r_n\rho_n t)^{\frac{\al}{2}} + h_1(n,s,t),
\label{a2}
\end{eqnarray}
where $h_1 = (r_n^2(t^2 + ||s||^2)) \frac{\al}{2} (\rho_n^2+
2r_n\rho_n t+\xi)^{\frac{\al}{2}-1},$ and $\xi \in (0,r_n^2(t^2 +
||s||^2)).$ Since $0 < \al \leq 2,$ and $(t,s) \in B(0,1),$ $0 \leq \xi
\leq r_n^2,$ and hence
\begin{eqnarray*}
 0 \leq  h_1(n,s,t) & = & (r_n^2(t^2 + ||s||^2))
\frac{\al}{2} (\rho_n^2+ 2r_n\rho_n t+\xi)^{\frac{\al}{2}-1}  \\
& \leq & (r_n^2(t^2 + ||s||^2))
\frac{\al}{2} (\rho_n^2+ 2r_n\rho_n t)^{\frac{\al}{2}-1}
  \leq   w_1(n),
\end{eqnarray*}
where
%
%
%
%
\begin{equation}
0  \leq w_1(n) := \frac{\al}{2} r_n^2 (\rho_n^2 -
2r_n\rho_n)^{\frac{\al}{2} - 1}   =
\frac{\al}{2} r_n^2 \rho_n^{\al-2}(1- \frac{2r_n}{\rho_n})^{\frac{\al}{2} - 1} \rar 0,
\label{w_1}
\end{equation}
since $r_n^{2}\rho_n^{\al-2} \rar 0,$ and $r_n/\rho_n \rar 0$ as $n \rar \infty.$
%
Again, from the Taylor's expansion applied to $(\rho_n^2+ 2r_n\rho_n t)^{\al/2}$ in
(\ref{a2}), we get
\begin{eqnarray}
(\rho_n^2+ 2r_n\rho_n t)^{\frac{\al}{2}} & = & \rho_n^{\al} + \al r_n
t\rho_n^{\al-1}+ \frac{\frac{\al}{2}(\frac{\al}{2}-1)}{2}(2r_nt \rho_n)^2(\rho_n^2+\zeta)^{\frac{\al}{2}-2}\nonumber\\
& = & \rho_n^{\al} + \al r_n t\rho_n^{\al-1} + h_2(n,t),\label{a3}
\end{eqnarray}
where $h_2(n,t) = \frac{\frac{\al}{2}(\frac{\al}{2}-1)}{2}(2r_nt
\rho_n)^2(\rho_n^2+\zeta)^{\frac{\al}{2}-2},$ and $\zeta \in (\min(0,2\rho_n
r_nt), \max(0,2 \rho_n r_nt)).$\\
Since $0 < \al \leq 2,$ and $-1 \leq t \leq 1,$ we get 
\begin{equation}
w_2(n) := \frac{\al(\al-2)}{2}
r_n^2\rho_n^{\al-2} \left(1 - 2 \frac{r_n}{\rho_n} \right)^{\frac{\al}{2} - 2} \leq
h_2(n,t) \leq 0.
\label{w_2}
\end{equation}
%
%
%
%
since $r_n^{2}\rho_n^{\al-2} \rar 0,$ and
$r_n/\rho_n \rar 0,$ it follows that $w_2(n)\rar 0$ as $n\rar \infty$.
%
From (\ref{a2})--(\ref{w_2}) we get
%
%
%
\begin{equation}
\rho_n^{\al} + 2\al r_n t \rho_n^{\al-1} + w_2 \leq ((\rho_n+r_nt)^2
+ (\|s\|r_n)^2 )^{\frac{\al}{2}} \leq \rho_n^{\al} + 2\al r_n t
\rho_n^{\al-1}+w_1.\label{a5}
\end{equation}
Using the above in (\ref{In}), we get
\begin{equation}
A_d r_n^d e^{-\lam w_1} G_n \leq I(\rho_n,r_n) \leq A_d r_n^d e^{-\lam w_2} G_n,
\label{bound_for_I}
\end{equation}
where
\begin{equation}
G_n = \int_{-1}^{1}\int_{\tiny{\|s\|^2 \leq(1-t^2),s\in
{\mathbb{R}}^{d-1}}} \exp\left(-\lam(\rho_n^{\al} + 2\al r_n t
\rho_n^{\al-1})\right) \ ds \ dt, \label{integral_in_bound_for_I}
\end{equation}
%
%
%
and $w_1$, $w_2$ as defined in (\ref{w_1}) and (\ref{w_2}) respectively, converge to $0$ as $n \rar \infty.$

If $\al >2,$ then $h_2(n,t) \geq 0,$ and we take $w_1,w_2$ to be the
sums of the upper and lower bounds of $h_1(n,s,t) + h_2(n,t)$ respectively.
We then obtain (\ref{bound_for_I}) with $w_2(n) = 0,$ and
%
\begin{eqnarray}
w_1(n) & = & \frac{\al}{2} r_n^2(\rho_n^2+2r_n\rho_n)^{\frac{\al}{2} - 1}
+ \frac{\al(\al-2)}{2} (r_n\rho_n)^2(\rho_n^2+2r_n\rho_n)^{\frac{\al}{2} - 2} \nonumber \\
 & = & \frac{\al}{2} r_n^2\rho_n^{\al-2}\left(1+2\frac{r_n}{\rho_n}\right)^{\frac{\al}{2}-1}
+ \frac{\al(\al-2)}{2} r_n^2\rho_n^{\al-2}\left(1+2\frac{r_n}{\rho_n}\right)^{\frac{\al}{2} - 2}
\label{pw_1}
\end{eqnarray}
%
%
%
%
which converges to zero by the conditions of the Lemma.

%
Now consider the integral in (\ref{integral_in_bound_for_I}). First make the change of variable
$u = t+1$ and then set $v = \lam \al r_n \rho_n^{\al-1}u$ to obtain
%
%
%
%
%
%
%
%
%
%
\begin{eqnarray}
G_n & = &  \theta_{d-1}  e^{-\lam \rho_n^{\al}}
\int_{-1}^{1}\exp(-\lam \al r_n \rho_n^{\al-1}t) (1-t^2)^{\frac{d-1}{2}} \ dt  \nonumber \\
 & = & \theta_{d-1} e^{-\lam (\rho_n^{\al}-\al r_n
\rho_n^{\al-1})} \int_{0}^{2}\exp(-\lam \al r_n \rho_n^{\al-1}
u)u^{(d-1)/2} (2-u)^{(d-1)/2} \ du \nonumber \\
& = &  \theta_{d-1} e^{-\lam (\rho_n^{\al}-\al r_n
\rho_n^{\al-1})}
(\lam \al r_n \rho_n^{\al-1})^{-\frac{d+1}{2}}
2^{\frac{d-1}{2}} K_n,
\label{a6}
\end{eqnarray}
%
where,
%
\begin{equation}
K_n = \int_{0}^{2\lam \al r_n \rho_n^{\al-1}}e^{-v} v^{\frac{d-1}{2}}
\left(1-\frac{v}{2\lam \al r_n
\rho_n^{\al-1}}\right)^{\frac{d-1}{2}} dv \leq
\Gamma((d+1)/2).\label{o1}
\end{equation}
We will show that as $r_n \rho_n^{\al-1} \rar \infty,$ the
integral in (\ref{a6}) converges to $\Gamma((d+1)/2)$ and also
estimate the error in this approximation. Write
\[ E_n := K_n - \Gamma((d+1)/2)) = A_n - B_n, \]
%
where,
\[A_n = \int_0^{2\lam \al r_n \rho_n^{\al-1}}  e^{-v}v^{(d-1)/2}\left[
\left(1-\frac{v}{2\lam \al r_n
\rho_n^{\al-1}}\right)^{(d-1)/2}-1\right]dv, \qquad \mbox{ and } \]
\[ B_n = \int_{2\lam \al r_n \rho_n^{\al-1}}^{\infty}e^{-v}v^{(d-1)/2} dv. \]
\[ \mid A_n\mid \leq \sup_{0\leq v\leq 2\lam \al r_n \rho_n^{\al-1}}
\left\{e^{-v/2}\mid 1 - \left( 1 - \frac{v}{2\lam \al r_n
\rho_n^{\al-1}} \right)^{(d-1)/2}\mid\right\} \int_0^{\infty}
e^{-v/2}v^{(d-1)/2}dv .\]

Since $(1-x)^a \geq 1 - Cx,$ $0 \leq x \leq 1$ with $C=1_{\{ 0<a \leq 1\}} + a 1_{\{a > 1\}},$ we
get,
\[0\leq 1 -  \left( 1 - \frac{v}{2\lam \al r_n \rho_n^{\al-1}} \right)^{(d-1)/2}
\leq\frac{C v}{2\lam \al r_n \rho_n^{\al-1}},\qquad
0\leq v \leq 2\lam \al r_n \rho_n^{\al-1}, \]
%
%
\[\Rar \qquad \mid A_n\mid \leq \frac{C}{2\lam \al r_n \rho_n^{\al-1}}
\; \sup_{0 \leq v < \infty} \left\{ v e^{-v/2} \right\}  \int_0^{\infty} e^{-v/2}v^{(d-1)/2}dv= \left(\frac{C'}{r_n
\rho_n^{\al-1}}\right), \]
where $C'$ is some constant. Further,
\[\mid B_n\mid \leq e^{-\lam \al r_n \rho_n^{\al-1}/2}\int_0^{\infty}e^{-v/2}v^{(d-1)/2} dv, \]
and hence decays exponentially fast in $r_n \rho_n^{\al-1}.$ Putting
the above two estimates in (\ref{a6}), we get
%
%
%
%
%
\begin{equation}
\mid E_n \mid \leq  \frac{C_1}{r_n \rho_n^{\al-1}} \rar 0, \qquad \mbox{ as } n \rar \infty.
\label{bound_for _E_n}
\end{equation}
%
%
%
%
The result now follows from (\ref{bound_for_I}), (\ref{a6}) and (\ref{bound_for _E_n}). \hfill $\Box$

We first prove Theorem~\ref{t1} for the number of isolated nodes
$W'_n(r_n)$ in the Poisson graph $G(P_n,r_n)$. Towards this end, we
first find an $r_n$ for which $E[W'_{n}(r_n)]$ converges. From the
Palm theory for Poisson processes (see (8.45),
Penrose (2003)), we get
\[E[W'_{n}(r_n)] = n \int_{R^d} \exp(-n I(x,r_n) f(x) dx.\]
Changing to Polar coordinates gives
\begin{equation}
E[W'_{n}(r_n)] = n \int_{0}^{\inft} \exp{(-nI(s,r_n))} f_R(s)ds,
\lab{pe2}
\end{equation}
where $f_R$ is defined in (\ref{e1}).
Let $\rho_n(t)^{\al} := \frac{t+a_n}{\lam},\: t \geq -a_n$ where
\begin{equation}
a_n := [\log \:n + (d/\al-1)\log_2n-\log(\Gamma{(d/\al)})].
\label{a_n}
\end{equation}
The idea is to make a change of variable $t = \rho_n^{-1}(s)$
such that $n f_R(\rho_n(t))\rho_n^{\prime}(t)$ converges and then
choose $r_n$ so that the first factor in (\ref{pe2}) also converges.
%
%
%
%
%
\begin{equation}
E[W_n'(r_n)] = \int_{-a_n}^{\infty} \exp(-n I(\rho_n(t),r_n)) g_n(t) dt, \label{pe2a}
\end{equation}
where
\begin{eqnarray}
g_n(t) & := & n f_R(\rho_n(t)) \rho_n'(t) = \frac{n \lam^{d/\al-1}}{\Gamma{(d/\al)}}
\left(\frac{t+a_n}{\lam}\right)^{\frac{d}{\al} -1}e^{-(t+a_n)}
 \nonumber\\
& = & \left( \frac{t+a_n}
{\log \:n}\right)^{\frac{d}{\al} -1} e^{- t} \nonumber\\
& = & \left( \frac{t+\log \:n +
(d/\al-1)\log_2n-\log(\Gamma{(d/\al)})}
{\log \:n}\right)^{\frac{d}{\al} - 1} e^{- t} \nonumber\\
%
%
& \rar & e^{- t}, \qquad \mbox{ as } n \rar \infty, \forall \;\; t
\in \mathbb{R}\ . \label{limit_g_n}
\end{eqnarray}
%
\begin{lem}
Suppose the sequence $\{r_n\}_{n\geq 1}$ satisfies
(\ref{rn}). Let $t \in \mathbb{R},$ and set $\rho_n(t)^{\al} =
\frac{t+a_n}{\lam} 1_{\{t \geq - a_n \}},$ where $a_n$ is as defined in (\ref{a_n}). Then
\begin{equation}
\lim_{n \rar \inft} nI(\rho_n,r_n) = C_d e^{\beta- t},\label{asy nI
}
\end{equation}
where $C_d$ is as defined in (\ref{e2a}).
\label{asymptotic radius}
\end{lem}
{\textbf{ Proof.}} It is easy to verify that for each fixed $t \in \mathbb{R}$, $\rho_n = \rho_n(t),r_n$ satisfy the conditions of Lemma~\ref{density over the ball} and so we have
%
\[
n I(\rho_n,r_n) \sim nA_d
\theta_{d-1}2^{\frac{d-1}{2}}\Gamma{(\frac{d+1}{2})} r_n^d
\exp(-\lam(\rho_n^{\al} - \al r_n\rho_n^{\al-1}))(\lam \al r_n
\rho_n^{\al-1})^{-\frac{d+1}{2}}. \]
Substituting for  $\lam\rho_n^{\al} = t +\log\:n +(\frac{d}{\al} -1)\log_2n
-\log(\Gamma(d/\al)),$ we get
\begin{equation}
nI(\rho_n,r_n)\sim \frac{nA_d\theta_{d-1} 2^{\frac{d-1}{2}}
\Gamma{(\frac{d+1}{2})}\Gamma{(d/\al)}e^{-t}}{n(\log\:n)^{\frac{d}{\al} -1}}r_n^{d}\exp(\lam\al
r_n\rho_n^{\al-1}))(\lam\al
r_n\rho_n^{\al-1})^{-\frac{d+1}{2}}.
\label{17.2}
\end{equation}
From (\ref{rn}), we can write
%
\begin{equation}
r_n = \frac{d-1}{\lam
\al}\frac{\log_2n}{(\lam^{-1}\log\:n)^{1-\frac{1}{\al}}} - \frac{d-1}{2\lam
\al}\frac{\log_3n}{(\lam^{-1}\log\:n)^{1-\frac{1}{\al}}} + \frac{\beta+
o(1)}{\lam \al(\lam^{-1}\log\:n)^{1-\frac{1}{\al}}},
\label{17.3}
\end{equation}
and hence
\begin{eqnarray}
\lam \al r_n\rho_n^{\al-1} & = &
\left(\frac{(d-1)\log_2n}{(\lam^{-1}\log\:n)^{1-\frac{1}{\al}}} -
\frac{d-1}{2}\frac{\log_3n}{(\lam^{-1}\log\:n)^{1-\frac{1}{\al}}} +
\frac{\beta+o(1)}{(\lam^{-1}\log\:n)^{1-\frac{1}{\al}}} \right)\nonumber\\
& \cdot & \left(\frac{1}{\lam}(t+\log\:n + (\frac{d}{\al}-1)\log_2n
-\log(\Gamma(d/\al)))\right)^{\frac{\al-1}{\al}}\nonumber\\
& = & \left((d-1)\log_2n -\frac{d-1}{2}\log_3n
+\beta+o(1)\right) \nonumber\\
&\cdot& \left(1+\frac{t}{\log\:n}+(\frac{d}{\al}-1)\frac{\log_2n}{\log\:n} -
\frac{\log(\Gamma(d/\al))}{\log\:n}\right)^{\frac{\al-1}{\al}} \label{17.35}\\
& = & (d-1)\log_2n - \frac{d-1}{2}\log_3n +\beta+o(1).\label{17.4}
\end{eqnarray}
Using (\ref{17.3}) and (\ref{17.4}) in  (\ref{17.2}), we
get
\begin{eqnarray*}
nI(\rho_n,r_n) & \sim & \frac{A_d\theta_{d-1} 2^{\frac{d-1}{2}}
\Gamma{(\frac{d+1}{2})}\Gamma{(d/\al)}e^{-t}}{(\log\:n)^{\frac{d}{\al} -1}}\\
& \cdot & \left(\frac{d-1}{\lam
\al}\frac{\log_2n}{(\lam^{-1}\log\:n)^{1-\frac{1}{\al}}} - \frac{d-1}{2\lam
\al}\frac{\log_3n}{(\lam^{-1}\log\:n)^{1-\frac{1}{\al}}} + \frac{\beta}{\lam
\al(\lam^{-1}\log\:n)^{1-\frac{1}{\al}}}\right)^d\\
& \cdot & \left(\frac{\exp((d-1)\log_2n - \frac{d-1}{2}\log_3n
+\beta))}{((d-1)\log_2n - \frac{d-1}{2}\log_3n
+\beta)^{\frac{d+1}{2}}}\right)\\
& = & A_d\theta_{d-1} 2^{\frac{d-1}{2}}
\Gamma{(\frac{d+1}{2})}e^{\beta-t}\left(\frac{d-1}{\lam
\al}\frac{\log_2n}{(\lam^{-1}\log\:n)^{1-\frac{1}{\al}}}
\left(1- \frac{\log_3n}{2\log_2n} + \frac{\beta}{(d-1)\log_2n}\right)\right)^d\\
& \cdot & \frac{(\log\:n)^{d-1}
(\log_2n)^{-\frac{d-1}{2}}}{(\log\:n)^{\frac{d}{\al} -1}}\Gamma(d/\al)
(d-1)^{-\frac{d+1}{2}}(\log_2n)^{-\frac{d+1}{2}}\\
& \sim & A_d\theta_{d-1} 2^{\frac{d-1}{2}}
\Gamma{(\frac{d+1}{2})}e^{\beta-t}\left(\frac{d-1}{\lam
\al}\frac{\log_2n}{(\lam^{-1}\log\:n)^{1- \frac{1}{\al}}}
\right)^d\\
& \cdot & (\log\:n)^{d-\frac{d}{\al} }
(\log_2n)^{-\frac{d-1}{2}}\Gamma(d/\al)
(d-1)^{-\frac{d+1}{2}}(\log_2n)^{-\frac{d+1}{2}} \;\; \rar \;\; C_de^{\beta-t}.
\end{eqnarray*}
\begin{lem}
There exists a constant $M$ depending on $\al, d$ and $\lam,$ such
that the following inequalities hold for all large enough $n$.
\begin{enumerate}
\item  Suppose $d/\al > 1,$ and $\lam r_n^{\al } - a_n \leq t \leq 0,$ or $d/\al < 1,$ and
$-\frac{\log\:n}{\log_2n} \leq t \leq 0,$ then $g_n(t) \leq M e^{-t}.$
\item For $d/\al < 1,$ and $\lam r_n^{\al } - a_n \leq t \leq -\frac{\log\:n}{\log_2n},$
$g_n(t) \leq M \left(\frac{\log_2n}{\log\:n}\right)^{d-\al}e^{-t}.$
\end{enumerate}
\label{lem_bounds_for_g_n}
\end{lem}

{\bf Proof. } Observe that for large $n$, $0.5 \log n \leq a_n \leq
2 \log n,$ and $\lam r_n^{\al} \geq \left( \frac{(d-1) \log_2
n}{2\al (\log n)^{1-1/\al}} \right)^{\al}.$

In the case when $d/\al > 1,$ and $\lam r_n^{\al } - a_n \leq t \leq
0,$
\[ g_n(t) \leq \left(\frac{0+a_n}{\log\:n}\right)^{\frac{d}{\al} -1}e^{-t} \leq 2^{\frac{d}{\al} -1}e^{-t}. \]
If $d/\al <1,$ and $-\frac{\log\:n}{\log_2n}\leq t \leq 0,$
\[ g_n(t) \leq \left(\frac{-\frac{\log\:n}{\log_2n}+ 0.5 \log n}{\log\:n}\right)
^{\frac{d}{\al} -1}e^{-t} \leq 4^{1 - \frac{d}{\al}} e^{-t}. \]
Finally, if $d/\al <1,$ and $\lam r_n^{\al}-a_n \leq t \leq
-\frac{\log\:n}{\log_2n},$
\[ g_n(t) \leq \left(\frac{\lam r_n^{\al}-a_n+a_n}{\log\:n}\right)^{\frac{d}{\al} -1}e^{-t}
 \leq \left(\frac{(d-1)\log_2n}{2 \al
\log\:n}\right)^{d-\al}e^{-t}. \qquad \qquad  \Box \]
We have the following proposition.
\begin{prop} Let the sequence $\{r_n\}_{n\geq 1}$ satisfy (\ref{rn}). Then
\eq{\lim_{n \rar \inft} E[W'_{n}] = \frac{e^{-\beta}}{C_d}\:,
\lab{pe3}}
where $C_d$ is as defined in (\ref{e2a}). \lab{prop1}
\end{prop}
{\textbf{ Proof.}} From Lemma~\ref{asymptotic radius} and (\ref{limit_g_n}), for each $t \in \mathbb{R}$,
we have
\begin{equation} \lim_{n \rar \inft}\exp(-n I(\rho_n(t),r_n))g_n(t) =
\exp(-C_d e^{\beta- t})e^{-t}.
\lab{e2b}
\end{equation}
%
%
Suppose we can find integrable bounds for $\exp(-nI(\rho_n(t),r_n))g_n(t)$ that hold
for all large $n.$ Then from (\ref{pe2a}), (\ref{e2b}) and the dominated
convergence theorem, we have
\begin{eqnarray*}
\lim_{n \rar \inft}E[W'_{n}(r_n)] & = & \lim_{n \rar \inft}
\int_{-a_n}^{\inft}\exp(-n I(\rho_n(t),r_n)) g_n(t) dt \\
& = & \int_{-\inft}^{\inft} \exp
\left(-C_de^{\beta-t}\right)e^{-t} dt = \frac{e^{-\beta}}{C_d}\int_{0}^{\inft}e^{-y} dy = \frac{e^{-\beta}}{C_d}. \label{e7}
\end{eqnarray*}
We find integrable bounds for $\exp(-nI(\rho_n(t),r_n))g_n(t)$, by dividing the range of $t$ into four parts.
\begin{enumerate}
\item First consider $t \geq 0.$ For large $n$ since $0.5 \log \: n < a_n < 2 \log \:n ,$ we have
\begin{equation}
g_n(t) \leq \left\{ \begin{array}{lc}
                              \left(\frac{(t+2\log \:n)}{\log \:n}\right)^{\frac{d}{\al}
                              -1}e^{-t}\ \leq  e^{-t} 2^{\frac{d}{\al}} \max(t,1)^{\frac{d}{\al} - 1}, & \qquad \frac{d}{\al} > 1, \\
                              \frac{e^{-t}}{2^{1- \frac{d}{\al}}} , & \qquad \frac{d}{\al} \leq 1.
                       \end{array} \right.
 \label{density bound}
\end{equation}
%
%
%
%
By the above bound on $g_n(t),$ it follows that
\begin{equation}
 \exp(-nI(\rho_n(t),r_n)) g_n(t) \leq g_n(t),  \lab{p4}
\end{equation}
is integrable over $(0,\infty).$
\item Now consider the range $-\frac{\log\:n}{\log_2n} \leq t \leq 0.$ As
$\lam \rho_n(t)^{\al} = t+a_n,$ from (\ref{17.35}) we get
\begin{eqnarray*}
\lam \al r_n\rho_n(t)^{\al-1} & = & ((d-1)\log_2n - \frac{d-1}{2}\log_3n +\beta + o(1)) \\
 & & \cdot \left(1 + \frac{\al -1}{\al} \left(\frac{t + (\frac{d}{\al}-1)\log_2n
-\log(\Gamma(d/\al))}{\log n} \right)(1 + \zeta_n(t))^{-\frac{1}{\al}}\right),
\end{eqnarray*}
where $\mid \zeta_n(t) \mid  \leq \mid t+(d/\al-1)\log_2n-\log(\Gamma{(d/\al)}) \mid (\log\:n)^{-1}.$
Since, $-\frac{\log\:n}{\log_2n}\leq t \leq 0,$ $\zeta_n(t) \rar 0,$ uniformly in
$t \in \left(-\frac{\log\:n}{\log_2n},0\right)$ as $n\rar \infty.$
Since $-1 \leq \frac{t \log_2 n}{\log n} \leq 0,$ in the above range
of $t$, we can find constants $c_1$ and $c_2$ such that for $n$ sufficiently large,
\begin{equation}
(d-1)\log_2n - \frac{d-1}{2}\log_3n - c_1 \leq \lam \al
r_n\rho_n(t)^{\al-1} \leq (d-1)\log_2n - \frac{d-1}{2}\log_3n + c_2.
\label{21.1}
\end{equation}
Hence for all sufficiently large $n$ we have
\begin{equation}
\exp(\lam \al r_n\rho_n^{\al-1}) \geq
\frac{(\log\:n)^{d-1}}{(\log_2n)^{\frac{d-1}{2}}}e^{-c_1}.\label{22.1}
\end{equation}
From Lemma~\ref{density over the ball},
\begin{eqnarray*}
\lefteqn{ nI(\rho_n,r_n)  \geq  nA_d\theta_{d-1}2^{\frac{d-1}{2}}
\left(\Gamma{(\frac{d+1}{2})}+E_n\right)r_n^de^{-\lam
w_1}\exp(-\lam(\rho_n^{\al}-\al r_n\rho_n^{\al-1}))(\lam
\al r_n\rho_n^{\al-1})^{-\frac{d+1}{2}} }\\
& = & nA_d\theta_{d-1}2^{\frac{d-1}{2}}
\left(\Gamma{(\frac{d+1}{2})}+E_n\right)r_n^d
\frac{\Gamma{(d/\al)}e^{-t}}{n(\log\:n)^{d/\al-1}}e^{-\lam
w_1}\exp(\lam \al r_n\rho_n^{\al-1})(\lam \al
r_n\rho_n^{\al-1})^{-\frac{d+1}{2}}.\label{22.2}
\end{eqnarray*}
Using (\ref{21.1}) and (\ref{22.1}) in above expression we
get
\begin{eqnarray*}
nI(\rho_n,r_n) & \geq & A_d\theta_{d-1}2^{\frac{d-1}{2}}
\left(\Gamma{(\frac{d+1}{2})}+E_n\right)e^{-\lam
w_1(n)}\left(\frac{(d-1)\log_2n}{\lam
\al (\lam^{-1}\log\:n)^{1-1/\al}}+o(1)\right)^d\nonumber\\
& \cdot & \frac{\Gamma{(d/\al)}e^{-t}}{(\log\:n)^{d/\al-1}}
\frac{(\log\:n)^{d-1}}{(\log_2n)^{\frac{d-1}{2}}}e^{-c_1}((d-1)\log_2n
- \frac{d-1}{2}\log_3n + c_2)^{-\frac{d+1}{2}}\nonumber \\
 & \geq & C \left( \Gamma{(\frac{d+1}{2})}+E_n \right) e^{-\lam
w_1(n)} e^{-t}.
\end{eqnarray*}
As in (\ref{21.1}), for $-\frac{\log\:n}{\log_2 n}\leq t\leq 0,$ it is easily verified that $r_n/\rho_n(t)$ and $r_n \rho_n(t)^{\al -2}$ converge uniformly to $0.$ It follows that $w_1(n)$ and $E_n$ converge uniformly to $0.$
Hence, we can find a constant $c' > 0,$ such that
\[ nI(\rho_n,r_n) \geq  c'e^{-t}. \]
From the above inequality and Lemma \ref{lem_bounds_for_g_n}(1),
there exists a constant $c$ such that for all sufficiently large
$n,$ we have
\begin{equation}
\exp(-nI(\rho_n(t),r_n))g_n(t) \leq c \exp(-c'e^{-t})e^{-t},\qquad
(-\frac{\log\:n}{\log_2n})\leq t \leq 0. \label{p4b}
\end{equation}
This upper bound is integrable over $ t \in(-\inft,0).$
\item Next, consider the range $\lam r_n^{\al}-a_n \leq t \leq
-\frac{\log\:n}{\log_2n}.$ From the first inequality we have $r_n \leq \rho_n(t),$ and hence
\begin{eqnarray}
I(\rho_n(t),r_n) & = & \int_{B(\rho_n(t)e,r_n)}A_d e^{-\lam \|x\|^{\al}}dx\\
& > & \int_{B(\rho_n(t)e,r_n), \|x\| \leq \rho_n(t)}A_d e^{-\lam\|x\|^{\al}}dx\\
& \geq & A_d e^{-\lam \rho_n(t)^{\al}}\mid B(\rho_n(t)e,r_n)\cap B(0,\rho_n(t))\mid, \label{aa1}
\end{eqnarray}
%
%
%
where $\mid \cdot \mid$ denotes the volume and $e =
(1,0,\ldots,0)\in \mathbb{R}^d.$ Inscribe a sphere of diameter $r_n$
inside $B(\rho_n(t)e,r_n)\cap B(0,\rho_n(t))$ (see Figure 1). Hence,
\begin{equation}
\mid B(\rho_n(t)e,r_n)\cap B(0,\rho_n(t))\mid \geq \frac{\theta_d
r_n^d}{2^d}.\label{a1}
\end{equation}
\begin{center}
  \includegraphics[scale=.50]{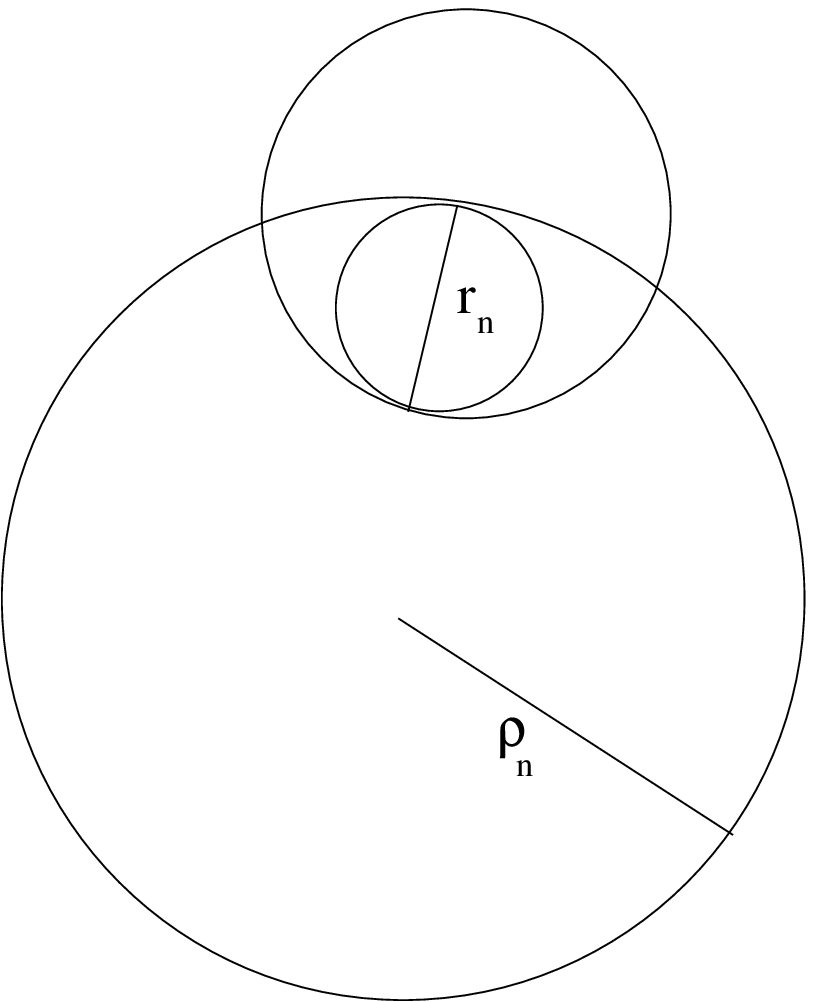}

{\bf Figure 1}
\end{center}
From (\ref{aa1}) and (\ref{a1}), we have
\begin{eqnarray}
I(\rho_n(t),r_n) &\geq & c'' e^{-\lam\rho_n(t)^{\al}} r_n^d\nonumber\\
& = & c'' \exp\left(-(t+\log\:n +(d/\al -1)\log_2 n
-\log(\Gamma{(d/\al)}))\right)\nonumber\\
& \cdot & \left(\frac{ (d-1)\log_2 n -
\frac{d-1}{2} \log_3 n +
\beta+o(1)}{ \lam \al (\lam^{-1}\log\:n)^{1-\frac{1}{\al}}}\right)^d \nonumber\\
& = & \frac{c''' e^{-t}}{n (\log\:n)^{\frac{d}{\al} - 1}}
\frac{ (\log_2 n)^d}{ (\log\:n)^{d - \frac{d}{\al}}}
\left( 1 - \frac{\log_3 n}{2 \log_2 n} + \frac{\beta + o(1)}{(d-1) \log_2 n}\right)^d \nonumber\\
& \geq & c^*n^{-1}(\log\:n)^{1-d}(\log_2n)^d e^{-t} = q_n e^{- t},\label{k1}
\end{eqnarray}
where $q_n = c^* (\log\:n)^{1-d}(\log_2n)^dn^{-1}$. From Lemma {\ref{lem_bounds_for_g_n}} and (\ref{k1}) we get,
%
\begin{eqnarray}
 \int_{\lam r_n^{\al}-a_n}^{-\frac{\log\:n}{\log_2n}}
\exp(-nI(\rho_n(t),r_n))g_n(t)dt
& \leq &
\left\{
  \begin{array}{ll}
    M \int_{\lam r_n^{\al}-a_n}^{-\frac{\log\:n}{\log_2n}}
\exp( - n q_n  e^{-t}) e^{-t} \: dt , & \hbox{$\frac{d}{\al}\geq 1$;} \\
 & \\
    M\left(\frac{\log_2n}{\log\:n}\right)^{d-\al} \int_{\lam r_n^{\al}-a_n}^{-\frac{\log\:n}{\log_2n}}
\exp( - n q_n  e^{-t}) e^{-t} \: dt , & \hbox{$\frac{d}{\al}< 1.$}
  \end{array}
\right. \nonumber  \\
 & & \nonumber  \\
& & \nonumber  \\
& \leq &
\left\{
  \begin{array}{ll}
    M \int_{\exp(\frac{
\log\:n}{\log_2n})}^{\exp(a_n - \lam r_n^{\al})} e^{-nq_n y} dy, & \hbox{$\frac{d}{\al}\geq 1$;} \\
 & \\
    M\left(\frac{\log_2n}{\log\:n}\right)^{d-\al}\int_{\exp(\frac{
\log\:n}{\log_2n})}^{\exp(a_n - \lam r_n^{\al})} e^{-nq_n y} dy, & \hbox{$\frac{d}{\al}< 1.$}
  \end{array}
\right. \nonumber  \\
& & \nonumber \\
& & \nonumber  \\
 & \leq &
\left\{
   \begin{array}{ll}
     \frac{M}{n q_n} e^{-nq_n e^{\frac{ \log\:n}{\log_2n}}}, & \hbox{$\frac{d}{\al}\geq 1$;}\nonumber\\
     \frac{M}{nq_n} \left(\frac{\log_2 n}{\log\:n}\right)^{d-\al}
e^{-nq_n e^{\frac{ \log\:n}{\log_2n}}}, & \hbox{$\frac{d}{\al}< 1.$}
   \end{array}
 \right. \\
\label{m1}
\end{eqnarray}
\begin{eqnarray}
\frac{M}{n q_n} e^{-nq_n e^{\frac{ \log\:n}{\log_2n}}} & = &
\frac{M}{nq_n}
\; \exp\left(-n^{1+\frac{1}{\log_2n}}q_n\right) \nonumber\\
%
%
& = & C \frac{(\log n)^{d-1}}{(\log_2n)^d}
\exp\left(-c^*n^{\frac{1}{\log_2n}}(\log\:n)^{1-d}(\log_2n)^{d}\right).
%
\label{m2}
\end{eqnarray}
Consider the exponent $c^* n^{\frac{1}{\log_2n}}(\log\:n)^{1-d}(\log_2n)^d.$ Taking logarithms,
we get
\[ \log(c^*)+\frac{\log\:n}{\log_2n}+ (1-d)\log_2n +d\log_3n
\geq \frac{\log\:n}{2\log_2n}. \]
Hence,
\begin{equation}
c^* n^{\frac{1}{\log_2n}}(\log\:n)^{1-d}(\log_2n)^d \geq
e^{\frac{\log\:n}{2\log_2n}}. \label{m5}
\end{equation}
Using (\ref{m5}) in (\ref{m2}), we get
\begin{equation}
\frac{M}{n q_n} e^{-nq_n e^{\frac{ \log\:n}{\log_2n}}} \leq
C \left(\frac{\log n}{\log_2n}\right)^{d-1} \frac{1}{\log_2n}
\exp\left(-e^{\frac{1}{2}\frac{\log\:n}{\log_2n}}\right)\rar 0,\label{m6}
\end{equation}
since the exponent is decaying exponentially fast in
$(\log\:n/\log_2 n).$
Using the inequality from (\ref{m6}) in (\ref{m1}) for the case $d/\al < 1,$ we get
\begin{equation}
\frac{M}{n q_n}\left(\frac{\log_2n}{\log\:n}\right)^{d-\al} e^{-nq_n
e^{\frac{ \log\:n}{\log_2n}}} \leq \frac{C(\log
n)^{\al - 1}}{(\log_2n)^{\al }}
\exp\left(-e^{\frac{\log\:n}{2\log_2n}}\right),\label{m7}
\end{equation}
which converges to $0,$ as $n \rar \infty,$ by the same argument as
above. From (\ref{m1}), (\ref{m6}) and (\ref{m7}) we have
\begin{equation}
 \int_{\lam r_n^{\al}-a_n}^{-\frac{\log\:n}{\log_2n}}
\exp(-nI(\rho_n(t),r_n))g_n(t)dt \rar 0.\label{e8}\\
\end{equation}
\item Finally, consider the case $-a_n\leq t \leq \lam r_n^{\al}-a_n.$ The
second inequality implies that $r_n \geq \rho_n(t).$ Hence for large
$n$ we have,
\begin{eqnarray}
nI(\rho_n(t),r_n) & = & n\int_{B(\rho_n(t) e,r_n)}A_d e^{-\lam \|x\|^{\al}}dx \nonumber \\
& > & n\int_{B(r_n e,r_n)}A_d e^{-\lam \|x\|^{\al}}dx \geq c_1 n e^{-\lam (2 r_n)^{\al}} r_n^d.
\label{lb_nI_final_case}
\end{eqnarray}
For large $n$ from (\ref{17.3}), we have
\begin{equation}
\frac{(d-1) \log_2 n}{2 \lam^{\frac{1}{\al}} \al (\log \: n)^{1 - \frac{1}{\al}}} \leq r_n \leq \frac{2 (d-1) \log_2 n}{\lam^{\frac{1}{\al}} \al (\log \: n)^{1 - \frac{1}{\al}}}.
\label{bounds_r_n}
\end{equation}
%
%
Fix $0 < \ep_1, \ep_2 < 1,$ such that $\ep = \ep_1 + \ep_2 < 1.$ Substituting from (\ref{bounds_r_n}) in (\ref{lb_nI_final_case}), we get, for large $n,$
\begin{eqnarray}
nI(\rho_n(t),r_n) & \geq & c_2 n e^{- c_3 \frac{(\log_2 n)^{\al}}{(\log \: n)^{\al - 1}}}
\frac{(\log_2 n)^{d}}{(\log \: n)^{d - \frac{d}{\al}}} \nonumber \\
 & \geq & c_2 n ^{1 - \ep_1}  e^{- c_3 \left( \frac{\log_2 n}{\log \: n} \right)^{\al} \log \: n}
\nonumber \\
 & = & c_2 n^{1 - \ep_1 - c_3 \left( \frac{\log_2 n}{\log \: n} \right)^{\al}} \geq c_2 n^{1 - \ep_1 -\ep_2}
= c_2 n^{1 - \ep}.
\label{final_bound_nI}
\end{eqnarray}
From (\ref{limit_g_n}), (\ref{final_bound_nI}) and the fact that for large $n$, $a_n < 2 \log \: n,$ we get
\begin{eqnarray}
\int_{-a_n}^{\lam r_n^{\al} - a_n} \exp(-nI(\rho_n(t),r_n)) g_n(t) \: dt
&  \leq &
\frac{e^{-c_2 n^{1-\ep}}}{(\log \: n)^{\frac{d}{\al} -1}}
\int_{-a_n}^{\lam r_n^{\al} - a_n} (t+a_n)^{\frac{d}{\al} -1} e^{-t} \: dt
\nonumber \\
 & \leq & \frac{e^{a_n}e^{-c_2 n^{1-\ep}}}{(\log \ n)^{\frac{d}{\al} -1}}
\int_0^{\infty} u^{\frac{d}{\al} -1} e^{-u} \: du \nonumber \\
 & \leq & \frac{ c n^2 e^{-c_2 n^{1-\ep}}}{(\log \ n)^{\frac{d}{\al} -1}}
\rar 0.
\label{e9}
\end{eqnarray}
\end{enumerate}
This completes the proof of Proposition~\ref{prop1}. \hfill $\Box$
%
%
\begin{thm}
Let $\al \in \mathbb{R}$ and let $r_n$ be as defined in (\ref{rn}). Then,
\[ W'_{0,n}(r_n) \stackrel{\cD}{\rar} Po(e^{-\beta}/C_d), \]
where $C_d$ is as defined in (\ref{e2a}) and and
$Po(e^{-\beta}/C_d)$ is the Poisson random variable with mean
$e^{-\beta}/C_d.$ \lab{limit_for_iso_vertices}
\end{thm}
{\textbf{ Proof.}} From Theorem~6.7, Penrose~(2003), we have
$d_{TV}(W'_{0,n}(r_n)$, $Po(E(W'_{0,n}(r_n)))$ is bounded by a constant times
$J_1(n)+J_2(n)$ where $J_1(n)$ and $J_2(n)$ are defined as follows.
\begin{equation}
J_1(n) = n^2\int_{\Re^d}\exp(-nI(x,r_n))f(x)dx \int_{B(x,3r_n)}
\exp(-nI(y,r_n))f(y)dy,
\label{eqn:def_J1}
\end{equation}
and
\begin{equation}
J_2(n) = n^2\int_{\Re^d}f(x)dx \int_{B(x,3r_n)\setminus B(x,r_n)
}\exp(-nI^{(2)}(x,y,r_n))f(y)dy,
\label{eqn:def_J2}
\end{equation}
where $I^{(2)}(x,y,r) = \int_{B(x,r)\cup B(y,r)}f(z)dz.$
Theorem~\ref{limit_for_iso_vertices} follows from Proposition~\ref{prop1} if we show that $J_i(n)
\rar 0$, as $n \rar \infty$, $i=1,2.$
We first analyze $J_1.$ Let $\rho_n(t),$ $g_n(t)$ be as defined in
Lemma~\ref{asymptotic radius} and (\ref{limit_g_n}) respectively.
\begin{eqnarray*}
J_1(n)
%
& = & n^2
\int_{-a_n}^{\infty}\exp(-nI(\rho_n(t),r_n))g_n(t) dt
\int_{B(\rho_n(t) e,3r_n)}\exp(-nI(y,r_n))f(y)dy\\
%
%
& = & J_{11}(n)+ J_{12}(n),
\end{eqnarray*}
where $J_{11}(n),$ and $ J_{12}(n)$ are defined as follows:
\begin{eqnarray*}
J_{11}(n) & = &
\int_{-a_n}^{-\frac{\log\:n}{\log_2n}}\exp(-nI(\rho_n(t),r_n))g_n(t)dt
\int_{B(\rho_n(t) e,3r_n)}\exp(-nI(y,r_n))n f(y)dy,\\
J_{12}(n) & = &
\int^{\infty}_{-\frac{\log\:n}{\log_2n}}\exp(-nI(\rho_n(t),r_n))g_n(t)dt
\int_{B(\rho_n(t) e,3r_n)}\exp(-nI(y,r_n)) n f(y)dy.
\end{eqnarray*}
First we will show that $J_{11}(n) \rar 0.$ From Proposition~\ref{prop1}, the inner integral in $J_{11}$,
\begin{eqnarray*}
\int_{B(\rho_n(t),3r_n)}\exp(-nI(y,r_n))n f(y)dy & \leq &
\int_{-a_n}^{\infty}\exp(-nI(\rho_n(t'),r_n))g_n(t')dt'  \\
 & & \;\; = E(W'_n(r_n)) \rar \frac{e^{-\beta}}{C_d}, \qquad \mbox{ as } n \rar
\infty.
\end{eqnarray*}
Thus, for any $\ep > 0,$ and all large $n,$ we have
\begin{equation}
J_{11}(n) \leq (1+ \ep) \frac{e^{-\beta}}{C_d}
\int_{-a_n}^{-\frac{\log\:n}{\log_2n}}\exp(-nI(\rho_n(t),r_n))g_n(t)dt.
\end{equation}
It follows from (\ref{e8}), (\ref{e9}) that $J_{11}(n)\rar 0.$
Next we will show that $J_{12}(n) \rar 0$ as $n \rar \infty.$
Define $B_n(t) = \{t': \rho_n(t) - 3 r_n \leq \rho_n(t') \leq \rho_n(t) + 3 r_n \}.$
The inner integral in $J_{12}(n),$
\begin{eqnarray*}
 \lefteqn{\int_{B(\rho_n(t)e,3r_n)}\exp(-nI(\rho_n(t'),r_n))g_n(t')dt'  } \\
%
& \leq &
 \left( 2 \sin^{-1}\left(\frac{3r_n}{\rho_n(t)}\right) \right)^{d-1} \int_{B_n(t)}\exp(-nI(\rho_n(t'),r_n))g_n(t')dt'\nonumber\\
& \leq &
\left( 2 \sin^{-1}\left(\frac{3r_n}{\rho_n(t)}\right) \right)^{d-1} \int_{-a_n}^{\infty}\exp(-nI(\rho_n(t'),r_n))g_n(t')dt'.\nonumber\\
& \leq &  2^{d-1}  (1+\ep) \frac{e^{-\beta}}{C_d}
\left( \sin^{-1}\left(\frac{3r_n}{\rho_n(t)}\right) \right)^{d-1} \leq C \left( \frac{\log _2 n}{\log \: n} \right)^{d-1},
\label{k2}
\end{eqnarray*}
since for all large $n,$ and $t \in (-\frac{\log\:n}{\log_2n}, \infty),$ we can find constants $c,c'$ and $\ep >0$ such that $0 \leq \frac{3r_n}{\rho_n(t)} \leq c \frac{\log _2 n}{\log \: n} \rar 0,$ and $\sin^{-1}(x) \leq c'x,$ for all $x \in [0, \ep].$
%
%
%
Thus the inner integral in $J_{12}$ converges uniformly to $0,$ as
$n \rar \infty.$ Hence $J_{12}$ converges to $0$ from the last
statement and the fact that the
expressions in (\ref{p4}), (\ref{p4b}) are integrable.

We now show that $J_2$ as defined in (\ref{eqn:def_J2}) converges to $0.$ Write
\begin{equation}
J_2(n) =J_{21}(n) + J_{22}(n) + J_{23}(n),
\label{j_2}
\end{equation}
where
\[ J_{2k}(n)  = n^2\int_{\Re^d}f(x)dx\int_{A_k(n)}\exp(-nI^{(2)}(x,y,r_n))f(y)dy,
\qquad k=1,2,3, \]
with $A_1(n) = \{ 2r_n \leq ||x-y|| \leq 3r_n\},\:$
$A_2(n) = \{r_n \leq ||x-y|| \leq 2r_n , \|x\| \leq \|y\|\},$ and
$A_3(n) = \{r_n \leq ||x-y|| \leq 2r_n , \|y\| \leq \|x\|\}.$
%
%
Since on $A_1(n)$, $I^{(2)}(x,y,r_n) = I(x,r_n) + I(y,r_n),$ we get,
\begin{eqnarray*}
J_{21}(n)
%
& = & n^2\int_{\mathbb{R}^d}\exp(-nI(x,r_n))f(x)dx\int_{\{ y: 2r_n \leq ||x-y|| \leq 3r_n\}}\exp(-nI(y,r_n))f(y)dy,\\
& \leq & n^2\int_{\mathbb{R}^d}\exp(-nI(x,r_n))f(x)dx
\int_{B(x,3r_n))}\exp(-nI(y,r_n))f(y)dy
= J_1(n),
\end{eqnarray*}
which has already been shown to converge to $0.$ Next we analyze  $J_{22}(n)$ as $n \rar \infty.$
The proof  for $J_{23}(n)$ is the same and so we omit it.
\begin{center}
  \includegraphics[scale=.30]{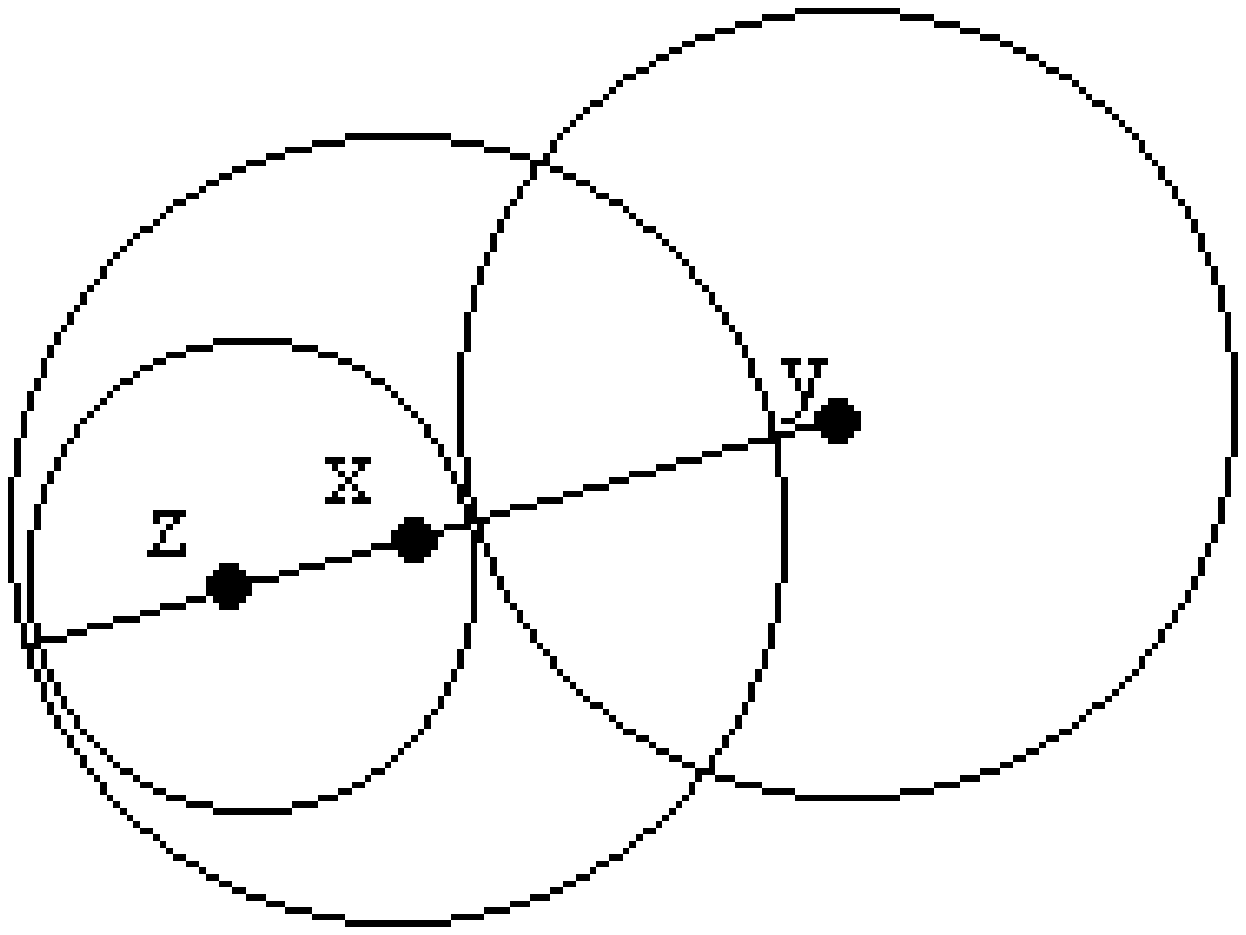}

{\bf Figure 2}
\end{center}
Let $B(z(x,y),\rho_1)$ be the ball with center $z = z(x,y)$ (see Figure 2) and radius $\rho_1 = \rho_1(x,y)\geq
\frac{r_n}{2}$ inscribed inside $B(x,r_n)\setminus B(y,r_n).$ Then
\begin{eqnarray*}
I^{(2)}(x,y,r_n) & \geq & I(z(x,y),\rho_1) + I(y,r_n)\\
& \geq & I(z(x,y),r_n/2) + I(y,r_n)\\
& \geq & I(x,r_n/2) + I(y,r_n),
\end{eqnarray*}
where the last inequality follows since $\|z \| < \|x\|.$ Thus,
\begin{eqnarray*}
J_{22}(n) & \leq & n^2 \int_{\mathbb{R}^d}
 \exp(-nI(x,r_n/2))f(x)dx\int_{A_2(n)}\exp(-nI(y,r_n))f(y)dy\\
& \leq & n^2\int_{\mathbb{R}^d}\exp(-nI(x,r_n/2))f(x)dx
\int_{B(x,3r_n))}\exp(-nI(y,r_n))f(y)dy\\
& = & J^*_{1}(n)+ J^*_{2}(n) + J^*_{3}(n),
\end{eqnarray*}
where
\begin{equation}
J_i^* = \int_{D_i}
\exp(-nI(\rho_n(t),r_n/2))g_n(t) \: dt
\int_{B(\rho_n(t),3r_n)}\exp(-nI(y,r_n))nf(y) \: dy, \qquad i=1,2,3,
\label{jstars}
\end{equation}

where $D_1 = [-a_n, -\frac{\log\:n}{\log_2n}),$ $D_2 = [ -\frac{\log\:n}{\log_2n}, 0)$ and
$D_3 = [0, \infty).$
The proof of $J_{i}^* \rar 0$, as $n \rar \infty$, for $i=1,3$ proceed exactly in the same manner as in
the case of $J_{11}$ and $J_{12}$ by replacing $r_n$ by $r_n/2$ while estimating
the outer integrals. In the case of $J_2^*$,  we proceed exactly as in the case of $J_{12}$
to obtain
\[ J_2^* \leq C \left( \frac{\log _2 n}{\log \: n} \right)^{d-1} \int_{D_2}
\exp(-nI(\rho_n(t),r_n/2))g_n(t) \:dt. \]
Estimating the integrand in the same way as in (\ref{p4b}), with $r_n$
replaced by $r_n/2$ and integrating, we get
\[ J_2^* \leq C' \left( \frac{\log _2 n}{\log \: n} \right)^{d-1} \left( \frac{(\log \: n)^{\frac{d-1}{2}}}{(\log_2 n)^{\frac{d-1}{4}}}
\right) \rar 0, \qquad \mbox{ as } n \rar \infty. \]
This completes the proof of Theorem~\ref{limit_for_iso_vertices}.\hfill $\Box$

{\textbf{Proof of Theorem~\ref{t1}.}} For each positive integer $n,$
set $m_1(n) =  n - n^{3/4}$ and $m_2(n) = n + n^{3/4}.$ Recall that the
Poisson sequence $N_n$ is assumed to be non decreasing.  Let $r_n$ be as in the
statement of the Theorem. It is easy to see that the proof of
Theorem~\ref{limit_for_iso_vertices} goes through for $m_i(n),$ that is,
\begin{equation}
W'_{m_i(n)}(r_n) \stackrel{\cD}{\rar} Po(e^{-\alpha}/C_d), \qquad i=1,2.
\lab{pe10}
\end{equation}
Let $\cP^-_n = \cP_{m_1(n)}$ and $\cP^+_n=\cP_{m_2(n)}.$ Let
$A^c$ denote the complement of set $A.$ Define events $H_n, A_n$ and
$B_n$ by
\begin{itemize}
\item $H_n = \{\cP^-_n\subseteq {\cal{X}}_n \subseteq
\cP^+_n\}.$
\item Let $A_n$ be the event that there exist a point $Y \in
\cP^+_n \backslash \cP^-_n$ such that $Y$ is isolated
in $G(\cP^-_n\cup\{Y\}, r_n).$
\item Let $B_n$ be the event that one or more points of
$\cP^+_n \backslash \cP^-_n$ lies within distance $r_n$
of a point $X$ of $\cP^-_n$ with degree zero in
$G(\cP^-_n,r_n).$
\end{itemize}
Then
\[ \{ W_{n}(r_n)\neq W'_{n}(r_n) \} \subseteq A_n \cup B_n \cup
F_n^c.\]
The proof is complete if we show that $P(A_n),P(B_n),P(F_n^c)$ all
converge to $0.$
\begin{eqnarray*}
P[H_n^c] & \leq & P[N_{m_1(n)} \geq n] + P[N_{m_2(n)} \leq n] \\
 & \leq & P[|N_{m_1(n)} - m_1(n)| \geq n^{3/4}] + P[|N_{m_2(n)} - m_2(n)| \geq n^{3/4}] \rar 0,
\end{eqnarray*}
by the Chebyshev's inequality.

Let $Y \sim f$ be a point independent of $\cP^-_n.$ Evidently,
\begin{eqnarray*}
P[A_n] & \leq & 2 n^{3/4} \: P[Y \mbox{ is isolated in }
G(\cP^-_n\cup \{Y\}, r_n)\}] \\
 & = &  2 n^{3/4} m_1(n)^{-1}E[W'_{m_1(n)}(r_n)] \rar 0, \qquad n \:
\rar \: \infty .
\end{eqnarray*}
By the Boole's inequality and the Palm theory,
\begin{eqnarray*}
P[B_n] & \leq & 2n^{3/4} \: P[\mbox{there is a isolated point of } G(\cP^-_n,r_n)
\mbox{ in } B(Y,r_n)] \\
 & \leq & 2 n^{7/4} \int_{\mathbb{R}^d}f(y)dy \cdot
\int_{B(y,r_n)} \exp(-m(n)I(x,r_n))f(x)dx.
\end{eqnarray*}
By interchanging the order of integration, we obtain
\begin{eqnarray}
P(B_n) & \leq & 2 n^{7/4} \int_{\mathbb{R}^d}I(x,r_n)\exp(-m(n)I(x,r_n))f(x)dx \nonumber\\
& = & 2 n^{3/4}\int_{-a_n}^{\infty}
I(\rho_n(t),r_n)\exp(-m(n)I(\rho_n(t),r_n))g_n(t)dt.
\label{bound_B_n}
\end{eqnarray}
From (\ref{density bound}) and (\ref{asy nI }), we get
\[2n^{3/4}I(\rho_n(t),r_n)\exp(-m(n)I(\rho_n(t),r_n))g_n(t) \leq C n^{-1/4}g_n(t)\rar 0. \]
Thus the integrand in (\ref{bound_B_n}) converges pointwise to $0$ as $n \rar \infty.$
Proceeding as in the proof of Proposition {\ref{prop1}}, using the
integrable bounds obtained in the proof of Proposition {\ref{prop1}},
for $\exp(-m(n) I(\rho_n(t),r_n)) g_n(t)$ and the bounds for
$I(\rho_n(t),r_n),$ and the dominated convergence theorem, we get
$P[B_n]  \rar 0.$ This completes the proof.\hfill $\Box$

{\textbf{Proof of Theorem~\ref{t2}.}} Let $r_n$ be as in the
statement of the Theorem. Then,
\[ \lim_{n \rar \infty} P[d_n \leq r_n] =
\lim_{n \rar \infty} P[W_{n}(r_n) = 0] = \exp(-e^{-\beta}/C_d).
 \]
\hfill $\Box$ \\
In order to prove strong law results for the LNND for graphs with
densities having compact support, one covers the support of the
density using an appropriate collection of concentric balls and then
shows summability of certain events involving the distribution of
the points of $\cX_n$ on these balls. The results then follow by an
application of the Borel-Cantelli Lemma. In case of densities having
unbounded support, the region to be covered changes with $n$ and
must be determined first. The following Lemma gives us the regions
of interest when the points in $\cX_n = \{X_1,X_2, \ldots ,X_n\}$, $
n\geq 1$ are distributed according to the probability density function
$f$ given by (\ref{d dim exponential density}) .

For any $c \in \mathbb{R}$, and large enough $n$, define
\eq{ R_n^{\al}(c) =  \frac{1}{\lam} \left(\log n +
\frac{c+d-\al}{\al}\log_2 n \right). \lab{def_R_n}}
For any set $A$, let $A^c$ denote its complement.
Let $U_n(c)$ be the event $\cX_n \subset B(0,R_n(c))$
and for any $c < 0$, $V_n(c)$ denote the event that at least one
point of $\cX_n$ lies in $B (0,R_n(0))\setminus B(0,R_n(c)).$ $a_n
\stackrel{>}{\sim} b_n$ implies that $a_n > c_n$ for some sequence
$c_n$ and $c_n \sim b_n.$ Further, $C,C_1,C_2,$ etc., will denote
constants whose values might change from place to place.
\begin{lem} Let the events $U_n$ and $V_n$, $n \geq 1,$ be as defined above. Then
\begin{enumerate}
\item $P[U_{n}^c(c) \mbox{ i.o. }] =
0,$ for any $c >\al,$ and
\item  $P[V_n^c(c) \mbox{ i.o. }] = 0,$ for any $c < 0$.
\end{enumerate}
The above results are also true with $\cX_n$ replaced by
$\cP_{\lam_n}$ provided $\lam_n \sim n.$ \lab{l2}
\end{lem}
Thus for almost all realizations of the sequence $\{\cX_n\}_{n \geq
1}$, all points of $\cX_{n}$ will lie within the ball $B(0,R_{n}(c))$ for any
$c> \al$ eventually, and for $c< 0$, there will be at least one
point of $\cX_n$ in $B(0,R_n(0))\setminus B(0,R_n(c))$ eventually.
%

%
{\bf Proof of Lemma~\ref{l2}. } From (\ref{e1}), write
$f_R(r)= A'_d e^{-\lam r^{\al}}r^{d-1}.$ Note that
\begin{equation}
\int_{\tilde{R}}^{\infty} f_R(r) \: dr \sim  A'_d (\lam \al)^{-1}  \tilde{R}^{d-\al} e^{-\lam \tilde{R}^{\al}},
\quad \mbox{ as }\: \tilde{R} \rar \infty. \label{eq3}
\end{equation}
Fix $a > 1,$ and define the subsequence $n_k = a^{k}.$ For large $k$, we have
\begin{eqnarray*}
P[\cup_{n= n_k}^{n_{k+1}}U_n^c(c)] & \leq &
P[\mbox{at least one vertex of ${\cal{X}}_{n_{k+1}}$ is in $B^c(0,R_{n_k}(c))$}] \\
& = & 1- (I(0,R_{n_k}(c)))^{n_{k+1}} = 1 - (1 - \int_{R_{n_k}(c)}^{\infty} f_R(r) \: dr)^{n_{k+1}} \\
& \leq & n_{k+1} \int_{R_{n_k}(c)}^{\infty} f_R(r) \: dr \sim
A'_d (\lam \al)^{-1} n_{k+1} R_{n_k}^{d-\al}(c) e^{-\lam R_{n_k}^{\al}(c)}  \\
& \leq & \frac{C}{k^{c / \al}}.
\end{eqnarray*}
Thus the above probability is summable for $c > \al,$ and the first
part of Lemma~\ref{l2} follows from the Borel-Cantelli Lemma.

Next, let $c< 0$ and take $n_k = a^k,$ for some $a > 1$. Note that
for all $n,m$ sufficiently large $R_n(c)$ are increasing and $R_n(c)
< R_m(0).$ Hence for $k$ sufficiently large, using (\ref{eq3}) and
the inequality $1-x \leq \exp(-x),$ we get
\begin{eqnarray}
P[\cup_{n= n_k}^{n_{k+1}} V_n^c(c)] & \leq &
P[\cX_{n_k} \cap (B(0,R_{n_k}(0))\setminus B(0,R_{n_{k+1}}(c))) = \emptyset] \nonumber \\
& = & \left(1 - \int_{R_{n_{k+1}}(c)}^{R_{n_k}(0)} A'_d e^{-\lam
r^{\al}}r^{d-\al} dr \right)^{n_k} \nonumber \\
& \leq & \exp\left( - n_k \int_{R_{n_{k+1}}(c)}^{R_{n_k}(0)} A'_d
e^{-\lam r^{\al}} r^{d-\al} dr \right) \nonumber \\
& \leq & \exp(- n_k c_1 A'_d(\lam \al)^{-1}(R_{n_{k+1}}^{d-\al}(c) e^{-\lam
R_{n_{k+1}}^{\al}(c)} -
R_{n_k}^{d-\al}(0)e^{-\lam R_{n_k}^{\al}(0)}))\nonumber \\
& \leq & e^{-c_2 k^{-c/\al}}
\end{eqnarray}
which is summable for all $c < 0.$ The second part of Lemma~\ref{l2}
now follows from the Borel-Cantelli Lemma. If $\cX_n$ is replaced be $\cP_{\lam_n},$ where $\lam_n \sim n,$ then
\begin{eqnarray*}
P[U_n^c(c)] & = & 1-\exp(-\lam_n(1-I(0,R_n(c))))\\
%
& \stackrel{<}{\sim} & \lam_n A'_d(\lam \al)^{-1}R_n^{d-\al}(c)\exp(-\lam
R_n^{\al}(c)) \\
& \sim & nA'_d(\lam \al)^{-1}R_n^{d-\al}(c)\exp(-\lam R_n^{\al}(c))),
\end{eqnarray*}
which is same as the $P[U_n^c(c)]$ in case of $\cX_n.$ Similarly,
one can show that $P[V_n^c(c)]$ has the same asymptotic behavior as
in the case of $\cX_n.$ Thus the results stated for $\cX_n$ also
hold for $\cP_{\lam_n}.$
\begin{prop} Let $t > d/\al \lam,$ and let
$r_n(t) = t (\lam^{-1}\log\:n)^{\frac{1}{\al} - 1} \log_2n.$ Then with
probability $1$, $d_n \leq r_n(t)$ for all large enough $n.$
\lab{prop2}
\end{prop}
{\bf Proof. } Let $c > \al$ and pick $u,t$ such that $(c+ \al (d
-1))/\al^2 \lam < u < t,$ and $\ep
> 0$ satisfying
\[ \ep + u < t. \]
From Lemma~\ref{l2}, $\cX_n \subset B(0,R_n(c))$ a.s. for all large
enough $n.$ For $m=1,2, \ldots,$ let $\nu(m) = a^m,$ for some $a >
1.$ Let $\ka_m$ (the covering number), be the minimum number of
balls of radius $r_{\nu(m+1)}(\ep)$ required to cover the ball
$B(0,R_{\nu(m+1)}(c)).$ For large $m,$ we have
\begin{eqnarray}
\ka_m & \leq & C_1 \frac{R_{\nu(m+1)}(c)^d}{r_{\nu(m+1)}^d(\ep)}\nonumber\\
& = & \frac{(\log(\nu(m+1)) + \frac{c+d-\al}{\al} \log_2(\nu(m+1)))^{d/\al}}
{\lam^d \ep^d(\log(\nu(m+1)))^{(d/\al-d)}(\log_2(\nu(m+1)))^d}\nonumber\\
& \leq & C_2
\left(\frac{m+1}{\log
(m+1)}\right)^d. \lab{p4a}
\end{eqnarray}
Consider the deterministic set $\{x_1^m,\ldots ,x_{\ka_m}^m\} \subset
B(0,R_{\nu(m+1)}(c)),$ such that
\[ B(0,R_{\nu(m+1)}(c)) \subset \cup_{i=1}^{\ka_m}
B(x_i^m,r_{\nu(m+1)}(\ep)).\]

Let $\al >1.$ Given $x \in \mathbb{R}^d,$ define $A_m(x)$ to be the
annulus $B(x,r_{\nu(m+1)}(u))\setminus B(x,r_{\nu(m+1)}(\ep)),$ and
let $F_m(x)$ be the event such that no vertex of
${\cal{X}}_{\nu(m)}$ lies in $A_m(x),$ i.e.
\begin{equation}
F_m(x) = \{{\cal{X}}_{\nu(m)}[A_m(x)]=0 \}, \label{fm}
\end{equation}
where ${\cal{X}}[B]$ denotes the
number of points of the finite set ${\cal{X}}$ that lie in $B.$
For any $x \in B(0,R_{\nu(m+1)}(c)),$ we have
\begin{eqnarray*}
P[X_i \in A_m(x)] & = & \int_{A_m(x)} f(y) \; dy \nonumber \\
& \geq & \int_{A_m(R_{\nu(m+1)}(c)e)} f(y) \; dy \nonumber \\
& = & I(R_{\nu(m+1)}(c), r_{\nu(m+1)}(u)) - I(R_{\nu(m+1)}(c),
r_{\nu(m+1)}(\ep)).
\end{eqnarray*}
Since $R_n(c), r_n$ satisfy the conditions of Lemma~\ref{density
over the ball}, we have for large $m,$
\begin{eqnarray*}
P[X_i \in A_m(x)] & \geq  & e^{-\lam
R_{\nu(m+1)}^{\al}(c)}(R_{\nu(m+1)}^{\al-1}(c))^{-\frac{d+1}{2}}\\
& \cdot & \left(c_1 e^{\lam \al
r_{\nu(m+1)}(u)R_{\nu(m+1)}^{\al-1}(c)}(r_{\nu(m+1)}(u))^{\frac{d-1}{2}}
- c_2 e^{\lam \al
r_{\nu(m+1)}(\ep)R_{\nu(m+1)}^{\al-1}(c)}(r_{\nu(m+1)}(\ep))^{\frac{d-1}{2}}\right)\\
& := & q_m.
\end{eqnarray*}
Substituting the values of
$R_{\nu(m+1)}(c)$ and $r_{\nu(m+1)}(\cdot)$ in $q_m$, we get for large $m$
\begin{equation}
q_m \leq (C(u) - C(\ep)) \frac{(\log (m+1))^{(d-1)/2}}{a^{m+1} (m+1)^{c/\al + d -\al
\lam u-1}}. \label{qn}
\end{equation}
Hence, for large $m$, we have
\eq{ P[F_m(x)] \leq  (1-q_m)^{\nu(m)} \leq
\exp(-\nu(m) q_m ) \leq \exp\left(-C\frac{(\log (m+1))^{(d-1)/2}}{
m^{\frac{c}{\al} + d -\al \lam u-1}}\right). \lab{p5}}
Set $G_m = \cup^{\ka_m}_{i=1}F_m(x_i^m).$ From (\ref{p4a}) and (\ref{p5}), we have for large $m,$
\begin{eqnarray}
P[G_m] & = & P[\cup^{\ka_m}_{i=1}F_m(x_i^m)] \leq  \sum^{\ka_m}_{i=1}P[F_m(x_i^m)]\nonumber\\
& \leq & C_2 \left(\frac{m+1}{\log
(m+1)}\right)^d \exp\left(-C\frac{(\log
(m+1))^{(d-1)/2}}{ (m+1)^{c/\al + d -\al \lam u-1}}\right),\nonumber\\
\end{eqnarray}
which is summable in $m$ since $u > \frac{c + \al(d - 1)}{\al^2
\lam}.$ By Borel-Cantelli Lemma, $G_m$ occurs only for finitely many
$m$ a.s.

Pick $n$, and take $m$ such that $a^m \leq n \leq a^{m+1}.$ If $d_n
\geq r_n(t)$, then there exists an $X \in \cX_n$ such that
$\cX_n[B(X,r_n(t))\setminus \{X\}] = 0.$ By Lemma~\ref{l2}, $X$ will be
in $B(0,R_{\nu(m+1)}(c))$ for all large enough $n,$ so there is some
$i \leq \ka_m$ such that $X \in B(x_i^m,r_{\nu(m+1)}(\ep)).$ So, if
$m$ is large enough,
\[ r_{\nu(m+1)}(\ep) + r_{\nu(m+1)}(u) \leq r_{\nu(m+1)}(t) \leq r_n(t).\]
So, $F_m(x_i)$ and hence $G_m$ occur. since $G_m$ occurs finitely
often a.s., $d_n \leq r_n(t)$ for all large $n,$ a.s. The result now
follows since $c > \al$ is arbitrary.

In the case when $\al \leq 1,$ cover the ball
$B(0,R_{\nu(m+1)}(c_1)),$ by the balls of radius $r_{\nu(m)}(\ep)$
and define the annulus $A_m(x)$ to be $B(x,r_{\nu(m)}(u))\setminus
B(x,r_{\nu(m)}(\ep)).$ Take $F_m(x) =
\{{\cal{X}}_{\nu(m+1)}[A_m(x)]=0 \}$ and proceed as in the case $\al
>1.$ This completes the proof of Proposition~\ref{prop2}.\hfill $\Box$

Now we derive a lower bound for $d_n$. Let $r_n(t) =
t\log_2n(\lam^{-1}\log\:n)^{1/\al-1}.$
\begin{prop} Let $t < (d-1)/\al \lam.$ Then with probability
$1$, $d_n \geq r_n(t),$ eventually. \lab{prop3} \end{prop}
{\bf Proof.} We prove the above proposition using the Poissonization
technique, which uses the following Lemma (see Lemma 1.4, Penrose
\cite{Penrose4}).
\begin{lem} Let $N(\lam)$ be Poisson a random variable with mean
$\lam.$ Then there exists a constant $c$ such that for all $\lam
> \lam_1,$
\[P[X > \lam+\lam^{3/4}/2]\leq c\exp(-\lam^{1/2}), \]
and
\[P[X < \lam-\lam^{3/4}/2]\leq c\exp(-\lam^{1/2}).\] \lab{l3}
\end{lem}

Enlarging the probability space, assume that for each $n$ there
exist Poisson variables $N(n)$ and $M(n)$ with means $n-n^{3/4}$ and
$2n^{3/4}$ respectively, independent of each other and of
$\{X_1,X_2,\ldots \}.$ Define the point processes
\[ \cP_n^- = \{X_1,X_2,\ldots,X_{N(n)} \}, \qquad \cP_n^+ = \{
X_1,X_2,\ldots , X_{N(n)+M(n)} \}. \]
Then, $\cP_n^-$ and $\cP_n^+$ are Poisson point processes on
$\mathbb{R}^d$ with intensity functions $(n-n^{3/4})f(\cdot)$ and
$(n+n^{3/4})f(\cdot)$ respectively. The point processes $\cP_n^-$,
$\cP_n^+$ and $\cX_n$ are coupled in such a way that $\cP_n^-
\subset \cP_n^+$. Thus, if $H_n = \{\cP_n^- \subset \cX_n \subset
\cP_n^+ \},$ then by the Borel-Cantelli Lemma and Lemma~\ref{l3},
$P[H_n^c \mbox{ i.o. }] = 0.$ Hence $\{\cP_n^- \subset \cX_n \subset
\cP_n^+ \}$ a.s. for all large enough $n.$

Pick constants $u,c,t, \ep$ such that $c < 0,$ $\ep > 0,$ $0 < t < u < (c+\al(d -1))/\al^2 \lam ,$
and $\ep + t < u.$

Consider the annulus $A_n(c) = B(0,R_n(0))\setminus B(0,R_n(c)),$ $c<0,$ where
$R_n(c)$ is as defined in (\ref{def_R_n}). For each $n,$ choose a
non-random set $\{ x_1^n,x_2^n,\ldots ,x_{\sg_n}^n \} \subset A_n(c),$
such that the balls $B(x_i^n,r_n(u)),$ $1 \leq i \leq \sg_n$ are
disjoint. The packing number $\sg_n$ is the maximum number of
disjoint balls $B(x,r_n(u))$, with $x \in A_n(c).$ For large $n,$ we have
\begin{eqnarray}
\sigma_n  & \geq &
c_1 \frac{R_n^d(0)-{R_n}^d(c)}{r_n^d(u)}\nonumber\\
& = & \frac{(\log\:n + \frac{d-\al}{\al}\log_2n)^{d/\al}
- (\log\:n + \frac{c+d-\al}{\al}\log_2n)^{d/\al}}{\lam^{d/\al}r_n^d(u)}\nonumber\\
& \geq & c_2 \left(\frac{\log \; n}{\log_2
n}\right)^{d-1}. \lab{p2}
\end{eqnarray}
By Lemma~\ref{l2}, there will be points in $A_n$ for all large
enough $n$, a.s.
Fix $a > 1$ and let $\nu(k) = a^k,$ $k=0,1,2, \ldots.$ 
Consider the sequence of sets
\begin{equation}
\left(\bigcup_{i=1}^{\sigma_{\nu(m)}}E_{m,i}\right)^c, \label{dn leq
rn}
\end{equation}
where,
\[E_{m,i} = \{{\cal{P}}^{-}_{\nu(m)}[B(x_i^{\nu(m)},r_{\nu(m)}(\ep))]= 1\}
\cap \{ {\cal{P}}^{+}_{\nu(m+1)}[B(x_i^{\nu(m)},r_{\nu(m)}(u))]= 1
\},
\]
where $i = 1,2,\ldots, \sg_{\nu(m)}, \:m=1,2,\ldots.$ From an
earlier argument $P[H_n^c]$ is summable and hence $H_n$ happens
eventually w.p.1. For any $n,$ let $m$ be such that $\nu(m) \leq n
\leq \nu(m+1).$ If $H_n$ and $E_{m,i}$ happen, then there is a point of $X \in
\cP^-_n \subset \cX_n$ such that $X \in B(x_i^{\nu(m)},r_{\nu(m)}(\ep))$ with no
other point of $\cP^+_n$ (and hence of ${\cal{X}}_n$) in $B(x_i^{\nu(m)}, r_{\nu(m)}(u)).$
This would imply that $d_n \geq r_{\nu(m)}(t) \geq r_n(t).$ Thus the
proof is complete if we show that

\[ \sum_{m=1}^{\infty} P\left[\left(\bigcup_{i=1}^{\sigma_{\nu(m)}}E_{m,i}\right)^c
\right] < \infty, \]
To this end, we first estimate $P[E_{m,i}].$

Set ${\cal{I}}_m = {\cal{P}}_{\nu(m+1)}^{+}\setminus
{\cal{P}}_{\nu(m)}^{-},$ and let $U_{m,i} =
B(x_i^{\nu(m)},r_{\nu(m)}(\ep)),$ and $V_{m,i} =
B(x_i^{\nu(m)},r_{\nu(m)}(u)) \setminus U_{m,i}.$ Then,
\[ E_{m,i} = \{{\cal{P}}^{-}_{\nu(m)}[U_{m,i}=1]\} \cap
\{{\cal{P}}^{-}_{\nu(m)}[V_{m,i}=0]\} \cap
\{{\cal{I}}_m[U_{m,i}=0]\} \cap \{{\cal{I}}_m[V_{m,i}=0]\}. \]
%
Let $\al(m) = \nu(m) - \nu(m)^{3/4}$ and $\be(m) = (\nu(m+1)) +
(\nu(m+1))^{3/4}.$ Note that each of the four events appearing in
the above equation are independent and that $\al(m) \sim \nu(m)$ and
$\be(m) \sim a \nu(m)$. Using this and the Lemma~\ref{density over
the ball}, we get, for all large enough $m,$
\begin{eqnarray}
P[E_{m,i}] & = & \al(m)\int_{U_{m,i}}f(y)\; dy \; \exp\left(
-\al(m)\int_{U_{m,i}}f(y) \; dy \right)  \exp \left(- \al(m)
\int_{V_{m,i}} f(y) \; dy \right)  \nonumber\\
& & \exp \left(-(\be (m)-\al(m)) \int_{U_{m,i}} f(y) \; dy \right)
\exp \left( - (\be(m)-\al(m)) \int_{V_{m,i}} f(y) \; dy \right) \nonumber\\
& = &  \al(m)\int_{U_{m,i}}f(y) \; dy \;
\exp \left( - \be(m) \int_{U_{m,i} \cup V_{m,i}} f(y)  \; dy \right)\nonumber\\
& = & \al(m) \; I(x_i^{\nu(m)},r_{\nu(m)}(\ep)) \; \exp
\left(- \; \be(m) \; I(x_i^{\nu(m)}, r_{\nu(m)}(u)) \right)\nonumber\\
& \geq & \al(m) \; I(R_{\nu(m)}(0),r_{\nu(m)}(\ep)) \; \exp
\left(- \; \be(m) \; I(R_{\nu(m)}(c), r_{\nu(m)}(u)) \right) \nonumber\\
& \geq &  C_1 \nu(m) r_{\nu(m)}^d(\ep)
e^{-\lam(R_{\nu(m)}^{\al}(0)-\al
r_{\nu(m)}(\ep)R_{\nu(m)}^{(\al-1)}(0))}(\lam \al
r_{\nu(m)}(\ep)R_{\nu(m)}^{(\al-1)}(0))
^{-\frac{d+1}{2}}\nonumber\\
& \cdot & \exp\left(-C_2 {\nu(m)} r_{\nu(m)}(u)^d
e^{-\lam(R_{\nu(m)}^{\al}(c)-{\al}
r_{\nu(m)}(\ep)R_{\nu(m)}^{\al-1}(c))}(\lam \al
r_{\nu(m)}(\ep)R_{\nu(m)}^{\al-1}(c))
^{-\frac{d+1}{2}}\right)\nonumber\\
& \sim & C_3 \frac{(\log_2(\nu(m)))^{\frac{d-1}{2}}}{(\log
(\nu(m)))^{d - 1 -\al \ep \lam}}
\exp\left(-C_4 \frac{(\log_2(\nu(m)))^{\frac{d-1}{2}}}{(\log(\nu(m)))^{d+c/\al-\al
u \lam-1}}\right)\nonumber\\
& \sim & C_3 \frac{(\log_2(\nu(m)))^{\frac{d-1}{2}}}{(\log
(\nu(m)))^{d - 1 -\al \ep \lam}}, \label{proba of en}
\end{eqnarray}
where the last relation follows since $u < \frac{\al d + c -\al
}{\al^2\lam}.$ The events $E_n(x_i^n),\:\:1\leq i\leq \sigma_n$ are
independent, so by (\ref{proba of en}), for large enough $m$,
\begin{eqnarray*}
P\left[\left(\bigcup_{i=1}^{\sigma_{\nu(m)}}E_{m,i}\right)^c\right]
& \leq & \prod_{i=1}^{\sigma_{\nu(m)}}\exp(-P[E_{m,i}]) \\
& \leq & \exp\left(-C_5 \sigma_{\nu(m)}\frac{(
\log_2(\nu(m)))^{\frac{d-1}{2}}}{(\log
(\nu(m)))^{d - 1 -\al \ep \lam}}\right) \\
& \leq & \exp\left(- C_6 \left(\frac{m}{\log\:m +\log
_2a}\right)^{d-1} \frac{(\log\:m +\log _2a)^{\frac{d-1}{2}}}{m^{d
- 1 -\al \ep \lam}}\right) \\
& = & \exp\left(-C_6 \frac{m^{\al \ep \lam}} {(\log\:m +\log
_2a)^{(d-1)/2}} \right),
\end{eqnarray*}
which is summable in $m.$

In case when $\al \leq 1.$ Define $U_{m,i} =
B(x_i^{\nu(m)},r_{\nu(m+1)}(\ep)),$ and $V_{m,i} =
B(x_i^{\nu(m)},r_{\nu(m+1)}(u)) \setminus U_{m,i}.$ Proceeding as
above we can show that
$P\left[\left(\bigcup_{i=1}^{\sigma_{\nu(m)}}E_{m,i}\right)^c\right]$
is summable. This gives $d_n \geq r_{\nu(m+1)}(t) \geq r_n(t).$
This completes the proof of Proposition~\ref{prop3}.\hfill $\Box$

{\textbf{Proof of Theorem~\ref{t3}.}} Immediate from
Proposition~\ref{prop2} and Proposition~\ref{prop3}.
\end{document}